\documentclass[12pt]{article}
\usepackage{amssymb}
\newcommand{\sect}[1]{\setcounter{equation}{0}\section{#1}}

\newcommand{\be}{\begin{equation}}
\newcommand{\ee}{\end{equation}}
\newcommand{\bea}{\begin{eqnarray}}
\newcommand{\eea}{\end{eqnarray}}
\newcommand{\beano}{\begin{eqnarray*}}
\newcommand{\eeano}{\end{eqnarray*}}

\newcommand{\vph}{\varphi}

\newcommand{\ca}{\mbox{$\cal{A}$}}
\newcommand{\cb}{\mbox{$\cal{B}$}}

\newcommand{\cf}{\mbox{${\cal F}$}}
\newcommand{{\cg}}{\mbox{$\cal{G}$}}

\newcommand{\cR}{\mbox{$\cal{R}$}}
\newcommand{\cs}{\mbox{$\cal{S}$}}

\newcommand{\prt}{\partial}

\newcommand{\wt}[1]{\widetilde{#1}}
\newcommand{\mb}[1]{\ \mbox{ #1 }\ }
\newcommand{\half}{\frac{1}{2}}

\newtheorem{prop}{Property}[section]

\newtheorem{defi}[prop]{Definition}
\newtheorem{theo}[prop]{Theorem}

\newtheorem{lem}[prop]{Lemma}
\newcommand{\prf}{\underline{Proof:}\ }
\newcommand{\finprf}{\null \hfill {\rule{5pt}{5pt}}\\[2.1ex]\indent}
\newcommand{\ie}{{\it i.e.}\ }
\newcommand{\CC}{\mbox{${\mathbb C}$}}
\newcommand{\ZZ}{\mbox{${\mathbb Z}$}}
\newcommand{\II}{\mbox{${\mathbb I}$}}

\begin{document}
\renewcommand{\thefootnote}{\fnsymbol{footnote}}
\newpage
\pagestyle{empty}
\setcounter{page}{0}

%%%%%%%%%%%%%%%%%%%%%%%%%%%%%%%%%%%%%%%%%%%%%%%%%
%%%%%%%%%%%%%%%%%%%% LOGO LAPTH - DEBUT  %%%%%%%%%%%%%%%
%%%%%%%%%%%%%%%%%%%%%%%%%%%%%%%%%%%%%%%%%%%%%%%%%
\newcommand{\LAP}{LAPTH}
\def\logo{{\bf {\huge LAPTH}}}

\centerline{\logo}

\vspace {.3cm}

\centerline{{\bf{\it\Large 
Laboratoire d'Annecy-le-Vieux de Physique Th\'eorique}}}

\centerline{\rule{12cm}{.42mm}}
%%%%%%%%%%%%%%%%%%%%%%%%%%%%%%%%%%%%%%%%%%%%%%%%
%%%%%%%%%%%%%%%%% LOGO LAPTH  - FIN %%%%%%%%%%%%%%
%%%%%%%%%%%%%%%%%%%%%%%%%%%%%%%%%%%%%%%%%%%%%%%%

\vspace{20mm}

\begin{center}
    
  {\LARGE  {\sffamily Vertex operators for boundary algebras }}\\[1cm]

{\large E. Ragoucy\footnote{ragoucy@lapp.in2p3.fr}}\\[.21cm] 
 Laboratoire de Physique Th{\'e}orique \LAP\footnote{UMR 5108 
    du CNRS, associ{\'e}e {\`a} l'Universit{\'e} de Savoie.}\\[.242cm]
    LAPP, BP 110, F-74941  Annecy-le-Vieux Cedex, France. 
\end{center}
\vfill\vfill

\begin{abstract}
We construct embeddings of boundary algebras $\cb$ into ZF algebras
$\ca$. 
Since it is known that these  algebras are the relevant ones for the study 
of quantum integrable systems (with boundaries for $\cb$ and without 
for $\ca$), this connection allows to make the link between different 
approaches of the systems with boundaries.
The construction uses the 
well-bred vertex operators built recently, and is classified 
by reflection matrices. It relies only on the existence of an 
$R$-matrix obeying a unitarity condition, and as such can be applied to 
any infinite dimensional quantum group.
\end{abstract}

\vfill
%\rightline{\tt mathQA/yymmnn}
\rightline{\LAP-860/01}
\rightline{July 01}

\newpage
\pagestyle{plain}
\setcounter{footnote}{0}
%%%%%%%%%%%%%%%%%%%%%%%%%%%%%%%%
%%%%%  FIN PAGE DE TITRE  %%%%%%
%%%%%%%%%%%%%%%%%%%%%%%%%%%%%%%%

%%%%%%%%%%%%%%%%%%%%%%%%%%%%%%%%
%%%%%  HEADINGS POUR DRAFT  %%%%%%
%%%%%%%%%%%%%%%%%%%%%%%%%%%%%%%%
\markright{\today\dotfill DRAFT\dotfill Operator vertex for 
boundary alg\dotfill }
%\pagestyle{myheadings}
%%%%%%%%%%%%%%%%%%%%%%%%%%%%%%%%
\sect{Introduction}
The problem of boundaries in integrable systems in the QISM framework 
was initiated by Cherednik \cite{Cher}. 
Basically, one can distinguish two approaches: the point of view of 
 Sklyanin \cite{Skly}, which relies 
on reflection matrices, and 
leads to the study boundary states (such as in \cite{GZ} or 
\cite{Cor} for instance); or the 
more algebraic approach of Mintchev et al \cite{Mint,bound}, where all the 
information is encoded in a boundary algebra $\cb$ \cite{bound}, 
and which allows to 
compute off-shell correlation functions and to study the integrals of 
motion \cite{BNLS}. 

On one hand, one starts with the bulk system and 
implement the boundary condition through a reflection matrix, while 
on the other hand, one has from the very beginning a boundary algebra 
$\cb$, 
which contains a reflection operator, and only after the 
specification of one of the several $\cb$-Fock space one gets a 
reflection matrix.

\null

In the present letter, we will take the second 
way of tackling the problem, and try to make more clear 
the connection with
 the first approach. For that purpose, we construct embeddings 
of $\cb$ into a Zamolodchikov-Faddeev (ZF) algebra, 
which is known to be the relevant algebra 
for the study of systems without boundary \cite{ZF}. 
These embeddings uses the 
well-bred vertex operators built in \cite{qOV}, and they are classified 
by the reflection matrices of the first approach, whence the link 
between the two points of view. 

Since one of the key point in our approach relies on the existence of 
so-called well-bred vertex operators, we present a generalization of this 
construction to the case of reflection operators.

The only assumption made for such  
constructions is the existence of an evaluated $R$-matrix (with 
spectral parameter) which obey 
the unitarity condition. It can thus be used for most of the integrable 
systems encountered in the literature.

\null 

The paper is organized as follows. In the  section \ref{s:def}, we 
introduce the different notions we will need. From these definitions, 
 we construct, in section \ref{s:prop}, a boundary algebra $\cb^{B}_{R}$
 from the deformed oscillator algebra $\ca_{R}$. 
Then, we consider, in section \ref{s:hier}, the hierarchy associated to 
$\cb^{B}_{R}$. Section \ref{s:ex} 
deals with one example: the nonlinear Schr\"odinger equation with 
boundary. Finally, 
we conclude in section \ref{s:concl}.

\sect{Definitions and notations\label{s:def}}
The starting point is an evaluated $R$-matrix, of size $N^2\times 
N^2$, with spectral parameter, and which obeys the Yang-Baxter equation 
and the unitarity condition:
\bea
&& R_{12}(k_1,k_2)R_{13}(k_1,k_3)R_{23}(k_2,k_3)=
R_{23}(k_2,k_3)R_{13}(k_1,k_3)R_{12}(k_1,k_2)\label{YBE}
\\
&& R_{12}(k_1,k_2)R_{21}(k_2,k_1)=\II\otimes \II \label{unitarity}
\eea
In the following, we will use the notation
\be
R_{12}=R_{12}(k_{1},k_{2}) 
\ee
\begin{defi}[ZF algebra $\ca_{R}$]\hfill\\
To the above $R$-matrix, one can associate a ZF algebra $\ca_{R}$, with generators 
$a_{i}(k)$ and $a^\dag_{i}(k)$ ($i=1,\ldots,N$) and exchange relations:
\bea
a_{1} a_{2} &=& R_{21} a_{2} a_{1} \label{AN-1}\\
a_{1}^\dag a_{2}^\dag &=& a_{2}^\dag a_{1}^\dag R_{21} \\
a_{1} a_{2}^\dag &=& a_{2}^\dag R_{12} a_{1} +\delta_{12}
\label{AN-3}
\eea
\end{defi}
We use the notations on auxiliary spaces
\bea
&&a_{1}=\sum_{i=1}^Na_{i}(k_{1})\, e_{i}\otimes\II
\mb{,} a_{2}=\sum_{i=1}^Na_{i}(k_{2})\, \II\otimes e_{i}\label{not1}\\
&&a^\dag_{1}=\sum_{i=1}^Na^\dag_{i}(k_{1})\, e^\dag_{i}\otimes\II
\mb{,} a^\dag_{2}=\sum_{i=1}^Na^\dag_{i}(k_{2})\, \II\otimes e^\dag_{i}\\
&& \delta_{12}=\delta(k_{1}-k_{2})\sum_{i=1}^N\, e_{i}\otimes e^\dag_{i} 
\mb{,} e^\dag_{i}=(0,\ldots,0,\stackrel{i}{1},0,\ldots,0)
\mb{,} e^\dag_{i}\cdot e_{j}=\delta_{ij}\ \ \label{not3}
\eea
where $\cdot$ stands for the scalar product of vectors.
\begin{defi}[Boundary algebra $\cb_{R}$] \hfill\\
\label{defB}
    To the same $R$-matrix, one can associate another algebra, 
    the boundary algebra $\cb_{R}$, 
    with generators $\wt{a}_{i}(k)$, $\wt{a}^\dag_{i}(k)$ and $b_{ij}(k)$
    ($i,j=1,\ldots,N$) and exchange relations:
    \bea
&&\wt{a}_{1}\wt{a}_{2} = R_{21} 
\wt{a}_{2}\wt{a}_{1} \label{BNl-1}\mb{;}
\wt{a}_{1}^\dag\wt{a}_{2}^\dag = 
\wt{a}_{2}^\dag\wt{a}_{1}^\dag R_{21}\\
&&\wt{a}_{1}\wt{a}_{2}^\dag = 
\wt{a}_{2}^\dag R_{12} \wt{a}_{1} +
\half{\delta_{12}}\delta+
\half{b_{12}}\label{BNl-3}\\
&&\wt{a}_{1}b_{2} = R_{21} 
b_{2}R'_{12} \wt{a}_{1}\mb{;}
b_{1}\wt{a}_{2}^\dag = 
\wt{a}_{2}^\dag R_{21}b_{1}R'_{21}\\
&& R_{12}\, b_{1}\, R'_{21}\, b_{2}\, =\,
b_{2}\, R'_{12}\, b_{1}\, \bar R_{21} \\
\label{rbrb}
&&b(k)b(-k) = \II 
\eea
\end{defi}
We have completed the notations (\ref{not1}-\ref{not3}) by:
\beano
&& R'_{12}=R_{12}(k_{1},-k_{2})\mb{;}R'_{21}=R_{21}(k_{2},-k_{1})\\
&& R'^{-1}_{12}=R_{12}(-k_{1},k_{2})\mb{;}
R'^{-1}_{21}=R_{21}(-k_{2},k_{1})\\
&& \bar R_{12}=R_{12}(-k_{1},-k_{2})\mb{;}
\bar R_{21}=R_{21}(-k_{2},-k_{1})\\
&&b_{12}=\delta(k_{1}+k_{2})\,
\sum_{i,j=1}^{N} b_{ij}(k_{1})\, e_{i}\otimes e_{j}^\dag\\
&& b_{1}(k)=\sum_{i,j=1}^{N} b_{ij}(k)\, E_{ij}\otimes \II_{N}\mb{;} 
b_{2}(k)=\sum_{i,j=1}^{N} b_{ij}(k)\, \II_{N}\otimes E_{ij}\, .
\eeano
Let us stress that
\be
R'^{-1}_{12}=\Big(R'_{12}\Big)^{-1}=\Big(R_{12}(k_{1},-k_{2})\Big)^{-1}=
R_{21}(-k_{2},k_{1})\neq R_{21}'
\ee
while $\bar R_{12}^{-1}=\bar R_{21}$.

The $\cb_{R}$ algebras have been introduced in \cite{bound}, where 
they were shown to 
play a fundamental role in the study of integrable systems with 
boundaries. They allow for instance the determination of off-shell 
correlation functions.

\null

Note that there is an automorphism on $\cb_{R}$ given by \cite{bound}:
\be
\rho\left\{\begin{array}{lcl}
\cb_{R} & \rightarrow & \cb_{R}\\[.23ex]
\wt{a}(k) & \mapsto & b(k)\wt{a}(-k)\\[.23ex]
\wt{a}^\dag(k) & \mapsto & \wt{a}^\dag(-k)b(-k)\\[.23ex]
b(k) & \mapsto & b(k)
\end{array}\right.\label{rho}
\ee
\begin{defi}[Reflection algebra $\cs_{R}$] \hfill\\
    The reflection algebra $\cs_{R}$ is the subalgebra of the boundary algebra, 
    with generators $b_{ij}(k)$
    ($i,j=1,\ldots,N$). It has exchange relations:
    \bea
&& R_{12}\, b_{1}\, R'_{21}\, b_{2}\, =\,
b_{2}\, R'_{12}\, b_{1}\, \bar R_{21} \\
&&b(k)b(-k) = \II 
\eea
\end{defi}
$\cs_R$ algebras enter into the class of $ABCD$-algebras introduced 
in \cite{MF}. They correspond, in the boundary algebra approach, 
to the symmetries of the underlying model (with boundary).

\begin{defi}[Well-bred vertex operator]\hfill\\
It has been shown in \cite{qOV}, that there exist in $\ca_{R}$  
a unique so-called well-bred vertex operator $T(k)=T^{ij}(k) E_{ij}$ such 
that
\bea
&& T(k_{\infty}) = \II+\sum_{n=1}^\infty\, \frac{(-1)^n}{(n-1)!} 
a^\dag_{{n}\ldots {1}}\, T^{(n)}_{\infty1\ldots n}a_{1\ldots n}\\
&& T_{1}a_{2} = R_{21}a_{2}T_{1} \label{defT}
\\
&& T_{1}a^\dag_{2} = a^\dag_{2}R_{12}T_{1}\\
&& R_{12}T_{1}T_{2} = T_{2}T_{1}R_{12} \label{rtt}
\eea
with 
\bea
a^\dag_{{n}\ldots {1}} &=& 
a^\dag_{\alpha_{n}}(k_{n})\ldots a^\dag_{\alpha_{1}}(k_{1})\\
a_{1\ldots n} &=& a_{\alpha_{1}}(k_{1})\ldots a_{\alpha_{n}}(k_{n}) \\
T^{(n)}_{\infty1\ldots n} &=& T^{(n)}_{\infty\alpha_1\ldots 
\alpha_n}(k_{\infty},k_1,\ldots,k_n)\in \left(\CC^{\otimes 
N^2}\right)^{\otimes (n+1)}\, (k_{\infty},k_1,\ldots,k_n)
\eea
In (\ref{defT}), there is an implicit summation on the indices
$\alpha_1,\ldots,\alpha_n=1,\ldots,N$ and an integration over the 
spectral parameters $\int\, dk_{1}\cdots dk_{n}$. The matrices
$T^{(n)}_{\infty1\ldots n}$ are built using only the evaluated $R$-matrix. 
For their exact expression, we refer to \cite{qOV}.
\end{defi}
\sect{Construction of $\cb_{R}$ from $\ca_{R}$\label{s:prop}}
\begin{theo}\label{bfroma}
Let $\ca_{R}$ be a ZF algebra, and $T$ its corresponding well-bred 
vertex operator. Let $B(k)$ be a $N\times N$ matrix such that
    \bea
&&R_{12}\, B_{1}\, R'_{21}\, B_{2}\, =\,
B_{2}\, R'_{12}\, B_{1}\, \bar R_{21}\label{RBRB}\\
&& B(k)B(-k)=\II_{N}
  \eea
Then, the following generators obey a boundary algebra $\cb_{R}^{B}$:
 \bea
\wt{a}(k) &=& \half\big(a(k)+b(k)a(-k)\big)\\
 \wt{a}^\dag(k) &=& \half\big( a^\dag(k)+a^\dag(-k)b(-k)\big)\\
 b(k) &=& T(k)B(k)T(-k)^{-1}
 \eea
 $B(k)$ is called the reflection matrix.
\end{theo}
\prf
We show the exchange relations by a direct calculation. 
As far as the exchange properties for $b$ are concerned, the proof 
follows the lines given in \cite{mol}, \ie a repetitive use of the 
relation (\ref{rtt}) in its various presentation, and also of
(\ref{RBRB}). We get
\be
R_{12}\, b_{1}\, R'_{21}\, b_{2}\, =\,
b_{2}\, R'_{12}\, b_{1}\, \bar R_{21}\label{eq:bb}
\ee
The unitarity condition $b(k)b(-k)=\II$ is obvious from the form of $b$.

We now compute the action of $b(k)$ on $a(k)$ and 
$a^\dag(k)$. For compactness, we write 
\be
a'=a(-k) \mb{;} a'^\dag=a^\dag(-k)\mb{;}T'=T(-k)
\mb{and} b'=b(-k)=b(k)^{-1}
\ee
\beano
a_{1}b_{2} &=& a_{1}T_{2} B_{2}T'^{-1}_{2} =
 R_{21}T_{2} a_{1}B_{2}T'^{-1}_{2}\\
&=& R_{21}T_{2} B_{2}a_{1}T'^{-1}_{2}=
 R_{21}T_{2} B_{2}T'^{-1}_{2}R'_{12}a_{1}
\eeano
Thus, we obtain
\be
a_{1}(k_{1})b_{2}(k_{2})=R_{21}(k_{2},k_{1})b_{2}(k_{2})
R_{12}(k_{1},-k_{2})a_{1}(k_{1})\mb{\ie}
a_{1}b_{2}=R_{21}b_{2}R'_{12}a_{1}
\label{eq:ab}
\ee
In the same way, we have
\beano
b_{1}a^\dag_{2} &=& T_{1} B_{1}T'^{-1}_{1}a^\dag_{2} =
 T_{1}B_{1} a^\dag_{2}T'^{-1}_{1}R_{21}'\\
&=& T_{1} a^\dag_{2}B_{1}T'^{-1}_{1}R_{21}'=
a^\dag_{2}R_{12}T_{1} B_{1}T'^{-1}_{1}R_{21}'
\eeano
that is
\be
b_{1}(k_{1})a^\dag_{2}(k_{2})=a^\dag_{2}(k_{2})R_{12}(k_{1},k_{2})
b_{1}(k_{1})R_{21}(k_{2},-k_{1})\mb{\ie}
b_{1}a^\dag_{2}=a^\dag_{2}R_{12}b_{1}R'_{21}
\label{eq:bad}
\ee

Thanks to the properties (\ref{eq:bb},\ref{eq:ab}), one computes
\beano
4\wt{a}_{1} \wt{a}_{2}  &=& a_{1}a_{2}+
a_{1}b_{2}a'_{2}+b_{1}a'_{1}a_{2}+b_{1}a'_{1}b_{2}a'_{2}\\
&=& R_{21} a_{2}a_{1}+R_{21}b_{2}R'_{12}a_{1}a'_{2}+
b_{1}R'_{21}a_{2}a'_{1}+(b_{1}R'_{21}b_{2}\bar R_{12})(a'_{1}a'_{2})\\
&=& R_{21} a_{2}a_{1}+R_{21}b_{2}a'_{2}a'_{1}+
R_{21}a_{2}b_{1}a'_{1}+ R_{21}b_{2}\, R'_{12}\, b_{1}\, \bar 
R_{21}a'_{2}a'_{1}\\
&=& R_{21} \Big(a_{2}a_{1}+b_{2}a'_{2}a'_{1}+
a_{2}b_{1}a'_{1}+ b_{2}a'_{2}\, b_{1}\, a'_{1}\Big)\\
&=& 4R_{21}\wt{a}_{2}\wt{a}_{1}
\eeano
The same calculation can be done for $\wt{a}^\dag_{1} \wt{a}^\dag_{2}$ using 
properties (\ref{eq:bb},\ref{eq:bad}).

We will use the unusual notation 
\be
\delta'_{12}=\delta(k_{1}+k_{2})\sum_{i=1}^N\, e_{i}\otimes 
e^\dag_{i}
\ee
Generally, the prime denotes a derivative when associated to the 
$\delta$ distribution, but fortunately the derivative of $\delta$ never 
occurs in the present article, so that there will be no confusion.

Looking now at $\wt{a}_{1} \wt{a}_{2}^\dag$, one gets:
\beano
4\wt{a}_{1} \wt{a}_{2}^\dag &=& a_{1}a_{2}^\dag+
a_{1}a'^\dag_{2} b'_{2}+b_{1}a'_{1}a_{2}^\dag+b_{1}a'_{1}a'^\dag_{2} b'_{2}\\
&=& a_{2}^\dag R_{12}a_{1}+\delta_{12}+
\big(a'^\dag_{2} R'_{12}a_{1} +\delta'_{12}\big)b'_{2}
+b_{1}\big(a_{2}^\dag a'_{1}+\delta_{12}'\big)
+b_{1}\big(a'^\dag_{2} \bar R_{12}a'_{1}+\delta_{12}\big)b'_{2}\\
&=& a_{2}^\dag R_{12}a_{1}+2\delta_{12}+2b_{12}+
a'^\dag_{2}b'_{2} R_{12}a_{1}+a_{2}^\dag R_{12}b_{1}a'_{1}
+a'^\dag_{2} R'_{21}b_{1}\bar R_{21}b'_{2}R'^{-1}_{21}a'_{1}\\
&=& 2\wt{a}_{2}^\dag R_{12}a_{1}+2\delta_{12}+2b_{12}+
a_{2}^\dag R_{12}b_{1}a'_{1}
+a'^\dag_{2} b'_{2}R_{12}b_{1}a'_{1}\\
&=& 4\wt{a}_{2}^\dag R_{12}\wt{a}_{1}+2\delta_{12}+2b_{12}
\eeano
Finally, from (\ref{eq:bb}), (\ref{eq:ab}) and (\ref{eq:bad}), one computes
the last equations, for instance:
\beano
2b_{1}\wt{a}_{2}^\dag &=& 
b_{1}\big({a}_{2}^\dag+{a'}_{2}^\dag b'_{2}\big)
={a}_{2}^\dag R_{12}b_{1}R'_{12}+{a'}_{2}^\dag R_{12}b_{1}\bar 
R_{21}b'_{2}\\
&=& {a}_{2}^\dag R_{12}b_{1}R'_{12}+{a'}_{2}^\dag 
b'_{2}R_{12}b_{1}R'_{21}=\wt{a}_{2}^\dag R_{12}b_{1}R'_{12}
\eeano
\finprf
\begin{lem}
    In $\cb_{R}^B$, the automorphism $\rho$ given in (\ref{rho}) is 
    the identity.
\end{lem}
\prf
Obvious, using $b(k)b(-k)=\II$.
\finprf
{\bf Remark:} Strictly speaking, $\cb_R^B$ corresponds to the coset 
(noted $\cb^\rho_R$) of the abstract boundary algebra $\cb_R$ (given by 
definition \ref{defB}) by the relation $\rho+id=0$, so that the construction of 
$\cb_R^B$ in theorem \ref{bfroma} defines 
inclusions of $\cb_R^\rho$ into $\ca_R$ (see also theorem \ref{ca/rho}).

\begin{prop}
 The Fock space $\cf_{R}$ of $\ca_{R}$ provides a Fock space representation 
 for $\cb_{R}^{B}$, defined by
\be
\wt{a}(k)\Omega=0\mb{and} b(k)\Omega=B(k)\Omega
 \ee
 \end{prop}
\prf
From the definition of the Fock space for $\ca_{R}$, one has 
$a(k)\Omega=0$, and thus $\wt{a}(k)\Omega=0$. 
The well-bred vertex operator 
$T(p)$ satisfies $T(p)\Omega=\Omega$, which implies 
$b(k)\Omega=B(k)\Omega$.
\finprf
{\bf Remark:} In \cite{bound,BNLS}, the boundary algebra $\cb_{R}$ has several 
Fock spaces, depending on the value of $b(k)$ on $\Omega$. In the 
present article, the algebra $\cb_{R}^B$ has only one Fock space, 
but $B$ is given within the construction of $\cb_{R}^B$, and there are as much 
$\cb_{R}^B$-algebra constructions in the present approach, as there are Fock spaces in 
the approach of \cite{Cher,Skly}.
\sect{Vertex operator construction}
As already mentioned, in \cite{qOV}, it has being shown that 
a well-bred vertex operator 
$T$ can be constructed as a series in $a$'s:
\be
T(k_{\infty}) = \II+\sum_{n=1}^\infty\, \frac{(-1)^n}{(n-1)!} 
a^\dag_{{n}\ldots {1}}\, T^{(n)}_{\infty1\ldots n}a_{1\ldots n}
\ee
with 
\bea
a^\dag_{{n}\ldots {1}} &=& 
a^\dag_{\alpha_{n}}(k_{n})\ldots a^\dag_{\alpha_{1}}(k_{1})\\
a_{1\ldots n} &=& a_{\alpha_{1}}(k_{1})\ldots a_{\alpha_{n}}(k_{n}) \\
T^{(n)}_{\infty1\ldots n} &=& T^{(n)}_{\infty\alpha_1\ldots 
\alpha_n}(k_{\infty},k_1,\ldots,k_n)\in \left(\CC^{\otimes 
N^2}\right)^{\otimes (n+1)}\, (k_{\infty},k_1,\ldots,k_n)
\eea
and an implicit summation on the indices
$\alpha_1,\ldots,\alpha_n=1,\ldots,N$ and an integration over the 
spectral parameters $\int\, dk_{1}\cdots dk_{n}$. 

The same expansion can be done for $T^{-1}$:
\be
T(k_{\infty})^{-1} = \II+\sum_{n=1}^\infty\, \frac{(-1)^n}{(n-1)!} 
a^\dag_{{n}\ldots {1}}\, T^{(n)\dag}_{\infty1\ldots n}a_{1\ldots n}
\label{Tinv}
\ee

For the exact expression of $T^{(n)}_{\infty1\ldots n}$ and
$T^{(n)\dag}_{\infty1\ldots n}$, we refer to 
\cite{qOV}. These matrices are constructed using only 
the evaluated $R$-matrix.

One can do the same construction for the $b$ operator:
\begin{prop}[Reflection operators as vertex operators]\hfill\\
In term of the $\ca_{R}$ generators, the reflection operators read
\beano
b_{0} &=& \II+\sum_{n=1}^\infty \frac{(-1)^n}{(n-1)!}
a^\dag_{{n}\ldots {1}}\, \beta^{(n)}_{01\ldots n}a_{1\ldots n}\\
\beta^{(n)}_{01\ldots n} &=& T^{(n)}_{01\ldots n}B_{0}+ 
B_{0}{T'}^{(n)\dag}_{01\ldots n}+(n-1)
\sum_{p=1}^{n-1}\left({n-2}\atop {p-1}\right) T^{(n-p)}_{0p+1\ldots n}
B_{0}{T'}^{(p)\dag}_{01\ldots p} 
\cR'_{p}
\eeano
where the prime ' indicates that one has to consider $-k_{0}$ instead 
of $k_{0}$ (as in definition \ref{defB}) and
\be
\cR'_{p}=R_{p+10}'^{-1}\cdots 
R_{n0}'^{-1}=\prod_{s=p+1}^{\longrightarrow\atop n}R_{0s}(-k_{0},k_{s})
\ee
\end{prop}
\prf
We start with $b=TBT'^{-1}$ and use the expansion of $T$. Then, from
$a_{1}T_{2}'^{-1}=T_{2}'^{-1}R_{12}'^{-1}a_{1}$, one gets
\be
b_{0}=B_{0}T_{0}'^{-1}+\sum_{n=1}^\infty \frac{(-1)^n}{(n-1)!}
a^\dag_{{n}\ldots {1}}\, T^{(n)}_{01\ldots n}B_{0}T_{0}'^{-1}
R_{10}'^{-1}\cdots R_{n0}'^{-1}a_{1\ldots n}
\ee
Finally, using the expansion (\ref{Tinv}) for $T_{0}'^{-1}$, 
and relabeling the 
auxiliary spaces, one obtains the result.
\finprf

Let us stress that the expansion is done in term of the $\ca_{R}$ 
generators $a$ and $a^\dag$, \underline{not} in term of 
$\wt{a}$ and $\wt{a}^\dag$, generators of $\cb_{R}$. It is possible 
that such an expansion would lead to a more simple expression for 
$\beta^{(n)}$.
\sect{Hierarchy for $\cb^{B}_{R}$\label{s:hier}}
\begin{prop}[Hierarchy for $\cb_{R}^B$]\hfill\\
    Let 
   \be
   H^{(n)}=\int_{-\infty}^\infty dk\, k^{n}\,  \wt{a}^\dag(k)\wt{a}(k)
   \ee
   Then:
   
   (i) $H^{(2n+1)}$ vanish identically in $\cb_{R}^B$.
   
   (ii) $[H^{(2n)},\wt{a}^\dag(k)]=k^{2n}\,\wt{a}^\dag(k)$ and 
   $[H^{(2n)},\wt{a}(k)]=-k^{2n}\,\wt{a}(k)$
   
   (iii) $\{H^{(2n)}\}_{n\in\ZZ_{+}}$ form a commuting flow for $\cb_{R}^B$, 
   called its hierarchy.
 \end{prop}
\prf
We first show $(i)$:
\beano
H^{(2n+1)}&=&\int_{-\infty}^\infty dk\, k^{2n+1}\, \wt{a}^\dag(k)\wt{a}(k)
=\int_{-\infty}^\infty dk\, (-k)^{2n+1}\, \wt{a}^\dag(-k)\wt{a}(-k)\\
&=& -\int_{-\infty}^\infty dk\, k^{2n+1}\, \wt{a}^\dag(k)b(k)b(-k)\wt{a}(k)
=-H^{(2n+1)}
\eeano
Now using the exchange relations of $\ca_{R}$, one gets
\beano
H^{(2n)}\wt{a}^\dag(k) &=& \int_{-\infty}^\infty dp\, p^{2n}\,  
\wt{a}_{1}^\dag(p)\wt{a}_{1}(p)\wt{a}^\dag_{2}(k)\\
&=& \int_{-\infty}^\infty dp\, p^{2n} \wt{a}_{1}^\dag(p)\Big(
\wt{a}^\dag_{2}(k)R_{12}(p,k)\wt{a}_{1}(p)+\half\delta(p-k)+\half 
b_{12}(p)\delta(p+k)\Big)\\
&=& 
\half\Big(k^{2n}\, \wt{a}^\dag(k)+(-k)^{2n}\, \wt{a}^\dag(-k)b(-k)\Big)\\
&&+
\int_{-\infty}^\infty dp\, p^{2n}\, \wt{a}^\dag_{2}(k)\wt{a}_{1}^\dag(p)
R_{21}(k,p)R_{12}(p,k)\wt{a}_{1}(p)\\
&=& k^{2n}\,\wt{a}^\dag(k) + \wt{a}^\dag(k)H^{(2n)}
\eeano
where in the last step we have used the automorphism $\rho=id$. The same 
computation leads to $[H^{(2n)},\wt{a}(k)]=-k^{2n}\, \wt{a}(k)$.

Note that starting with $H^{(2n+1)}$ and performing the above 
calculation leads to e.g. 
$[H^{(2n+1)},\wt{a}(k)]=0$, which is compatible with $(i)$.

Finally, using $(ii)$, one computes
\beano
{[H^{(2n)},H^{(2m)}]} &=& \int_{-\infty}^\infty dp\, p^{2n}\Big(
\wt{a}^\dag(p)\,[\wt{a}(p),H^{(2m)}]+\wt{a}(p)\,[\wt{a}^\dag(p),H^{(2m)}]
\Big)\\
&=& \int_{-\infty}^\infty dp\, p^{2n}\Big(
p^{2m}\wt{a}^\dag(p)\wt{a}(p)-p^{2m}\wt{a}^\dag(p)\wt{a}(p)\Big)\\
&=& 0
\eeano
\finprf
{\bf Remark:} On the Fock space $\cf_{R}$, and considering the 
states $\wt{\vert k}>=\wt{a}^\dag(k)\Omega$, one has
\be
H^{(2n)}\Omega=0\ \Rightarrow\  
H^{(2n)}\wt{\vert k}>=[H^{(2n)},\wt{a}^\dag(k)]\Omega=
k^{2n} \wt{\vert k}>
\ee
which shows that $H^{(2n)}\neq0$.

\null

{\bf Remark 2:} The hierarchy defines integrable systems with boundary 
defined by $B(k)$. In the framework we have adopted, the definition of 
the boundary is given by the data of the 
reflection matrix $B(k)$, as it is presented in \cite{Cher}, 
\underline{but} the boundary algebra is naturally recovered here, 
contrarily to \cite{Cher}, where it is lacking for the calculation of 
off-shell correlation functions. On the other hand, in 
\cite{bound,Mint,BNLS}, the boundary 
algebra is the basic data (whence the possibility of computation of 
correlation functions), \underline{but} the data of the boundary 
condition (\ie the reflection matrix) is given with the choice of a Fock space 
$\cf_{R}^{B}$. Thus, the present framework can be viewed as a bridge 
between the approaches \cite{Cher} and \cite{bound,BNLS}.

This remark is confirmed in the following theorem:

\begin{theo}\label{ca/rho}
Let $B$ be a reflection matrix of $\ca_R$, and $b=TBT'^{-1}$ the 
corresponding reflection operator. Let $\rho_B$ be defined by
\be
\rho_B(a)=b\,a'\mb{and}\rho_B(a^\dag)=a'^\dag\,b'
\ee
Then:

(i) $\rho_B$ is an automorphism of $\ca_R$

(ii) $\cb^B_R$ is the coset of $\ca_R$ by the ideal $Ker(\rho_B-id)$
\end{theo}
\prf
We prove $(i)$ by direct calculation, using the results 
(\ref{eq:ab},\ref{eq:bad}) and the exchange relations 
(\ref{BNl-1},\ref{rbrb}). Let 
$\alpha=\rho_B(a)$ and $\alpha^\dag=\rho_B(a^\dag)$:
\beano
\alpha_1\alpha_2 &=& b_1a'_1b_2a'_2=b_1R'_{21}b_2\bar R_{12}a'_1a'_2
=R_{21}b_2R'_{12}b_1\bar R_{21}a'_2a'_1\\
&=& R_{21}b_2R'_{12}R'^{-1}_{12}a'_2b_1a'_1=R_{21}\alpha_2\alpha_1
\eeano
The same calculation can be done with $\alpha^\dag_1\alpha^\dag_2$. For 
the last relation, we have:
\beano
\alpha_1\alpha^\dag_2 &=& b_1a'_1a'^\dag_2b'_2
=b_1(a'^\dag_2\bar R_{12}a'_1+\delta_{12})b'_2\\
&=& a'^\dag_2 R'_{12}b_1\bar R_{21}\bar R_{12}\bar 
R_{21}b'_2 R'^{-1}_{21}a'_1+\delta_{12}
= a'^\dag_2 b'_2 R_{12} b_1 R'_{21} R'^{-1}_{21} a'_1+\delta_{12}\\
&=& \alpha^\dag_2 R_{12} \alpha_1 +\delta_{12}
\eeano
The proof for $(ii)$ follows from the observation that 
$\wt{a}=\half\big(a+\rho_B(a)\big)$ and 
$\wt{a}^\dag=\half\big(a^\dag+\rho_B(a^\dag)\big)$. This means
\be
a\equiv \wt{a}\ \ \mbox{mod}[Ker(\rho_{B}-id)]
\mb{and}
a\equiv \wt{a}\ \ \mbox{mod}[Ker(\rho_{B}-id)]
\ee
 so that making the 
coset by $Ker(\rho_B-id)$ we recover $\cb_R^B$. 
\finprf 
\null

We present now a property which was already proved in \cite{BNLS} for the 
case of additive $R$-matrix $R_{12}(k_{1}-k_{2})$, but the reasoning is valid 
in full generality:
\begin{prop}[Integrals of motions of the hierarchy]\hfill\\
 The reflection algebra $\cs_{R}^B$ generates integrals of motion for 
 the $\cb_{R}^B$-hierarchy.
\end{prop}
\prf
By direct computation using the exchange relations (\ref{BNl-3}), one 
shows that $[H^{(2n)},b(k)]=0$.
\finprf
Still following the lines given in \cite{BNLS}, one gets:
\begin{prop}[Spontaneous symmetry breaking]\hfill\\
In the Fock space representation, there is a spontaneous symmetry 
breaking of the symmetry algebra through
\be
b(k)\Omega=B(k)\Omega\label{bo=Bo}
\ee
\end{prop}
\prf
From (\ref{bo=Bo}), one knows all
the operators $b_{ij}(k)$ which have non-vanishing value on $\Omega$.
Since the $\cs_{R}^{B}$-algebra constitutes the symmetry algebra of our
problem, we are exactly faced with a mechanism of spontaneous
symmetry breaking for our reflection algebra.
\finprf
\sect{Example: the nonlinear Schr\"odinger equation with 
boundary\label{s:ex}}
It has already been shown \cite{Mint,BNLS} that all the informations on the 
hierarchy associated to the  nonlinear Schr\"odinger  equation in 1+1 
dimensions with boundary (BNLS) can be reconstructed starting from a 
boundary algebra $\cb_{R}$, 
where $R$ is the $R$-matrix of the Yangian $Y(N)$ based on $gl(N)$
\be
R(k)=\frac{1}{k+ig}\left(k\,\II_{N}\otimes \II_{N} 
+ig\, P_{12}\right)\mb{,} P_{12}=\sum_{i,j=1}^N E_{ij}\otimes E_{ji}
\ee
This $R$-matrix obey an additive Yang-Baxter equation
\be
R_{12}(k_1-k_2)R_{13}(k_1-k_3)R_{23}(k_2-k_3)=
R_{23}(k_2-k_3)R_{13}(k_1-k_3)R_{12}(k_1-k_2)
\ee
and one shows, using $P^2=\II$, that $R_{12}(k)R_{21}(-k)=\II$. Thus, 
the properties stated above apply.

In fact, it is well-known that
 the canonical field $\Phi$ obeying the (quantum) NLS:
\[
\Big(i\prt_t+\prt^2_x\Big)\Phi(x,t) = 2g\, 
:\Phi(x,t)\bar\Phi(x,t)\Phi(x,t):
\mb{with} \Phi(x,t)=\left(\begin{array}{c} \vph_{1}(x,t) \\ \vdots \\ 
\vph_{n}(x,t) \end{array}\right)
\]
has an Hamiltonian which is exactly 
$H^{(2)}$, and that the reflection algebra $\cs_{R}$ is a symmetry of the 
hierarchy. 

To make the contact with the present point of view, one has to specify 
the Fock space $\cf_{R}^{B}$, which amount to fix a boundary matrix $B$. Then, 
one can construct the boundary algebra $\cb_R^B$, 
which will have the same Fock space. The data of $B$ also completely 
determine the boundary condition for the physical field $\Phi$ 
(see \cite{BNLS}), as it should in the approach of \cite{Skly}.
The phenomenon of spontaneous symmetry 
breaking in BNLS was also studied in \cite{BNLS}.

\sect{Conclusion\label{s:concl}}
We have shown that one can embed the boundary algebra $\cb_{R}$
into the ZF algebra $\ca_{R}$, and that there are as much embedding as 
there are reflection matrices. Such embeddings allow to link the 
approach of Mintchev et al \cite{bound,Mint}, who introduced the boundary 
algebras, to the original work of Cherednik \cite{Cher} and lately 
Sklyanin \cite{Skly}, who studied the problem of factorized S-matrices 
in models with boundaries. In particular, the results presented here 
allow to reconstruct the boundary algebra from the ZF algebra and the 
data of a reflection matrix.

The construction relies only on an $R$-matrix with spectral parameter 
which satisfies the Yang-Baxter equation and a unitarity condition. As a 
consequence, it is applicable to most of the infinite dimensional 
quantum groups, and in particular to Yangians,
and to centerless affine or elliptic quantum groups.

Taking as an example the $R$-matrix of $Y(N)$, the Yangian based on $gl(N)$, 
we recover by this construction the nonlinear Schr\"odinger equation 
with boundary and its symmetry. It is thus very natural to believe that the 
other integrable systems known in the literature can be treated with 
the present approach.

As an extension to our approach, the generic problem of (elliptic) 
quantum algebras with central extension should also be treated. 

%\newpage

\end{document}